\pdfoutput=1

\documentclass[final, reqno, 10pt]{amsart}

 \usepackage{amsmath,amsthm,amsfonts}
\usepackage{latexsym,mathabx}
\usepackage{amssymb}
\usepackage{bbold}
\usepackage{enumerate}

\newtheorem{theorem}{Theorem}[section]
\newtheorem{lemma}[theorem]{Lemma}
\newtheorem{prop}[theorem]{Proposition}
\newtheorem{corollary}[theorem]{Corollary}
\newtheorem{cor}[theorem]{Corollary}
\theoremstyle{remark}

\newtheorem{definition}{Definition}[section]
\newcommand{\set}{\mathbb}

\newcommand{\les}{\lesssim}
\renewcommand{\frak}{\mathfrak}

\newcommand{\mc}{\mathcal}
\newcommand{\be}{\begin{equation}}
\newcommand{\ee}{\end{equation}}
\newcommand{\bee}{\begin{align}}
\newcommand{\eee}{\end{align}}
\newcommand{\ba}{\begin{array}}

\newcommand{\ea}{\end{array}}
\newcommand{\bpm}{\begin{pmatrix}}
\newcommand{\epm}{\end{pmatrix}}
\newcommand{\lb}{\label}

\def\calH{\mc{H}}

\def\wtil{\widetilde}

\DeclareMathOperator{\Imim}{Im}

\newcommand{\oast}{\circledast}

\newcommand{\ov}{\overline}
\newcommand{\dd}{{\,}{d}}

\renewcommand{\Im}{\Imim}
\newcommand{\R}{\mathbb R}

\newcommand{\Z}{\mathbb Z}
\def\Sph{\mathbb{S}}

\newcommand{\B}{\mathcal B}

\newcommand{\one}{{\mathbb{1}}}
\newcommand{\EQ}[1]{\begin{equation}  \begin{split} #1 \end{split} \end{equation} }
\def\embed{\hookrightarrow}
\def\const{\mathrm{const}\cdot}
\def\eps{\varepsilon}
\def\what{\widehat}
\def\mes{\mc{M}}
\def\calB{\mc{B}}
\def\calF{\mc{F}}
\def\calS{\mc{S}}

\def\nn{\nonumber}
\def\trip{|\!|\!|}
\def\intR{\int_{-\infty}^{\infty}}
\def\wtild{\widetilde}
\def\cB{\calB}

\renewcommand{\epsilon}{\varepsilon}

\numberwithin{equation}{section}

\begin{document}

\title[Structure formulas for wave operators, scaling invariant norm]{Structure formulas for wave operators under a   small scaling invariant  condition}

\author{M.\ Beceanu}
\thanks{The first author thanks the University of Chicago for its hospitality during the summers of 2015 and 2016.}
\address{Department of Mathematics, University at Albany, State University of New York, 1400 Washington Avenue, Albany, NY 12222 }
\email{mbeceanu@albany.edu}

\author{W.\ Schlag}
\thanks{The second author was partially supported by NSF grant DMS-1500696
during the preparation of this work.}
\address{Department of Mathematics, The University of Chicago, 5734 South University Avenue, 
Chicago, IL 60637}
\email{schlag@math.uchicago.edu}

\begin{abstract}
We continue our work on  the structure formula for the intertwining wave operators $W_{\pm}$ associated with $H=-\Delta+V$ in $\R^{3}$, cf.~\cite{BeSc}.  
We consider small potentials relative to a scaling invariant norm. 
\end{abstract}

\maketitle

\section{Introduction}
\lb{sec:intro}

In a recent paper~\cite{BeSc} we obtained a structure formula for the intertwining wave operators  
 for $H=-\Delta+V$ in three dimensions. We imposed the following condition on the potential $V$. Define $B^{\beta}$, $\beta\ge0$, as   $L^{2}$ functions with  
 \EQ{\label{eq:Bsp}
\|f\|_{B^{\beta}}:= \| \one_{[|x|\le 1]} f\|_{2} + \sum_{j=0}^{\infty} 2^{j\beta} \big\| \one_{[2^{j}\le |x|\le 2^{j+1}]} f \big\|_{2} <\infty
}
Then for $V$ real-valued, $V\in B^{\beta}(\R^{3})$, $\beta\ge\frac12$,  the wave operators 
\[
W_{\pm} = \lim_{t\to\pm \infty} e^{itH} e^{-itH_{0}}
\]
exist in the strong $L^{2}$ sense, with $H_{0}=-\Delta$. These operators satisfy for continuous, bounded $f$ on the line, $f(H) W_{\pm} = W_{\pm} f(H_0)$, and $P_c=W_\pm W_\pm^*$,  where 
$P_c$ is the projection onto the absolutely continuous spectral subspace of $H$ in $L^{2}$. There is no singular continuous spectrum (asymptotic completeness).    Yajima~\cite{yajima0,yajima}, established the $L^{p}$ boundedness of the wave operators assuming that zero energy is neither an eigenvalue nor a resonance.  

In \cite{BeSc} we proved the following theorem. By $B^{1+}$ we mean $B^{\beta}$ for some $\beta>1$. By $ \mc M_y $ we mean the Borel measures on Euclidean space. 

\begin{theorem}[\cite{BeSc}]
\lb{thm:struct}
Let $V \in B^{1+}$ be real-valued and assume that $H = -\Delta+V$ admits no eigenfunction or resonance at zero energy. Then there exists  $g(x, y, \omega) \in L^1_\omega \mc M_y L^\infty_x$, i.e., 
\EQ{\nn
  \int_{\Sph^2} \|g(x, dy, \omega)\|_{\mc M_y L^{\infty}_x}  \dd \omega < \infty
} 
such that for any $f \in L^2$ one has the representation formula
\be\begin{aligned}\nn
(W_{+} f)(x) &= f(x) + \int_{\Sph^2} \int_{\R^3} g(x, dy, \omega) f(S_\omega x - y)   \dd \omega.
\end{aligned}\ee
where $S_{\omega}x=x-2(x\cdot\omega)\omega$ is a reflection. A similar result holds for $W_{-}$. 
\end{theorem}

As an application, suppose $X$ is any Banach space of measurable functions on $\R^3$ which is invariant under translations and reflections, and in which Schwartz functions are dense. Assume   that  $\|\one_H  f\|_X\le A\|f\|_X$ for all half spaces $H\subset\R^3$ and $f\in X$ with some uniform constant~$A$. Then 
\EQ{\label{eq:HW}
\| W_+ f\|_X \le A\, C(V) \|f\|_X \qquad\forall\; f\in X
}
where $C(V)$ is a constant depending on $V$ alone.  In particular, this recovers Yajima's $L^p$ boundedness of the wave operators.  Furthermore, \cite{BeSc} obtains quantitative estimates on the norm in~\eqref{eq:HW} as well as on $C(V)$. These bounds blow up as $\beta\to 1$ with $V\in B^\beta$. Nevertheless, we remark that one can obtain Theorem~\ref{thm:struct}, albeit  without quantitative control, under the condition $V\in B^1$ although the details are not worked out in~\cite{BeSc}. 

The goal here is to seek a scaling invariant condition on $V$ under which a structure formula~\eqref{eqn1.8} can be obtained. The natural scaling of the Schr\"{o}dinger operator $H=-\Delta+V$ is 
$V\to \lambda^2 V(\lambda x)$, $\lambda>0$ in any dimension. In the framework of the $B$-spaces above~\eqref{eq:Bsp} the critical norm relative to this scaling is $\dot B^{\frac12}(\R^3)$ where 
 \EQ{\label{eq:Bsp dot}
\|f\|_{\dot B^{\beta}}:= \sum_{j=-\infty}^{\infty} 2^{j\beta} \big\| \one_{[2^{j}\le |x|\le 2^{j+1}]} f \big\|_{2} 
}
This norm is invariant under  the aforementioned scaling provided $\lambda$ is a power of~$2$. 
However, it is currently unclear whether Theorem~\ref{thm:struct} might hold for potentials $V\in \dot B^{\frac12}$. It is possible that the threshold $\beta=1$ could be optimal for the $B^\beta$ spaces. 
It is natural to investigate the scaling invariant class for several reasons: (i) it is an optimal scenario, corresponding to the $|x|^{-2}$ decay rate which balances the Laplacian (ii) it arises in widely-studied energy critical
nolinear equations such as the $u^5$ wave equation in $\R^3$:
\[
u_{tt} -\Delta u - u^5=0
\]
which admits explicit the $1$-parameter stationary solutions $u(t,x)= \lambda^{\frac12} W(\lambda x)=:W_\lambda(x)$, $\lambda>0$, $W(x) = (1+|x|^2/3)^{-\frac12}$.  In the radial class, these are the only stationary nonzero solutions of finite energy. Linearizing about $W_\lambda$ yields the family of Schr\"{o}dinger operators $H_\lambda=-\Delta - \lambda^2 W^4(\lambda x)$. It is therefore desirable to work with a condition on the potential that is uniform in $\lambda>0$. 

This paper presents such a norm, but currently we only consider small potentials in this norm. To formulate it, we recall some notation.   

\begin{definition}
\label{def:trip}
For any Schwartz function $V$ we define $\trip V\trip=\|L_V\|_{L^1_{t, \omega}}$, where
\[
L_V(t,\omega)  =  \int_0^\infty \what{V}(-\tau\omega) e^{\frac{i}{2}t\tau}\, \tau\, d\tau
\]
is as above. For any Schwartz function $v$ in $\R^3$
\EQ{\lb{eq:Bnorm2}
\|v\|_B := \sup_{\Pi} \int_{-\infty}^\infty \trip\delta_{\Pi(t)}\, v(x)\trip \, dt 
}
where $\Pi$ is a $2$-dimensional plane through the origin, and $\Pi(t)=\Pi+t\vec N$, $\vec N$ being the unit norm to $\Pi$. 
\end{definition} 

Clearly,  $\trip v\trip \le \|v\|_{B}$. We will show below that $\| v\|_B<\infty$ is finite for Schwartz functions. We will do this by dominating $\| v\|_B$ 
by a stronger norm which is also scaling invariant and more explicit, cf.\ Lemma~\ref{lem:Bnorm dom}. This more explicit norm involves half of a derivative on $2$-planes. So it is not
a pure decay condition on the potential. 
The   key analytical arguments in this paper are based on  the   
precise norm as defined in~\eqref{eq:Bnorm2}.  The main result of this paper is the following one.

\begin{theorem} 
\lb{thm:main} 
There exists $c_0>0$ with the following property: for any real-valued $V$ with 
 $\| V\|_B + \| V\|_{\dot B^{\frac12}}\le c_0$,  there exists  $g(x, y, \omega) \in L^1_\omega \mc M_y L^\infty_x$ 
 with 
 \EQ{\lb{eq:gfin}
  \int_{\Sph^2} \|g(x, dy, \omega)\|_{\mc M_y L^{\infty}_x}  \dd \omega \les c_0
} 
such that for any $f \in L^2$ one has the representation formula
\be\begin{aligned}\lb{eqn1.8}
(W_{+} f)(x) &= f(x) + \int_{\Sph^2} \int_{\R^3} g(x, dy, \omega) f(S_\omega x - y)   \dd \omega.
\end{aligned}\ee
where $S_{\omega}x=x-2(x\cdot\omega)\omega$ is a reflection. A similar result holds for $W_{-}$. 
\end{theorem}

 The spectral properties of $H=-\Delta+V$ are irrelevant under the smallness assumption.  For $L^{1}\to L^{\infty}$ dispersive estimates of the Schr\"{o}dinger evolution $e^{itH}P_{c}$ with a scaling-invariant condition on~$V$   in $\R^{3}$ without any smallness assumption,  see~\cite{becgol}.  It is not clear to the authors if there might be other scaling-invariant norms which are better suited for structure theorems for the wave operators. This question is particularly relevant with respect to large potentials and the Wiener formalism that is instrumental for  Theorem~\ref{thm:struct}.  

\section{The wave operators and their expansion}
\label{sec:wav}

We now recall the formalism of the wave operator $W_{+}$ going back to Kato~\cite{kato}.   %
First, by \cite[Lemma~2.2]{BeSc}, $\dot B^{\frac12}\embed L^{\frac32,1}$  (using Lorentz space notation).  If $V\in L^{\frac32,1}$, then  the wave operators 
$W_{\pm}$ exist, and are isometries from $L^{2}$ onto the range of $P_{c}$, the projection onto the continuous spectrum of $H$. Moreover, if $V$ is small in $ L^{\frac32,1}$, then $H$ has no eigenvalues and no zero energy resonance, and the spectrum is purely absolutely continuous. In other words, $W_{\pm}$ are unitary operators.  Moreover, for any $f\in L^{2}(\R^{3})$ the integral
\be\begin{aligned}\label{eq:WDuh}
W_+ f &= f + i \int_0^{\infty} e^{itH} V e^{-itH_0} f \dd t 
\end{aligned}\ee
converges in the strong sense.   See for example Section~4 of~\cite{BeSc} for more details. 

Expanding  \eqref{eq:WDuh} iteratively by means of the Duhamel formula one has 
\EQ{ \lb{Dyson}
W_+ f &= f + W_{1+} f + \ldots + W_{n+} f + \ldots,\\
  W_{1+} f &= i \int_{t>0} e^{-i t \Delta} V e^{i t \Delta} f \dd t,\ \ldots \\
W_{n+} f &= i^n \int_{t>s_1>\ldots>s_{n-1}>0} e^{-i(t-s_1)\Delta} V e^{-i(s_1-s_2) \Delta} V \ldots \\
  &e^{-i s_{n-1} \Delta} V e^{it\Delta} f \dd t \dd s_1 \ldots \dd s_{n-1}
}
 for $f\in L^2$.  For small potentials one can actually sum this series, which will give Theorem~\ref{thm:main}. 
In addition to the operators $W_{n+}$,  we shall work with their regularized version, 
\be\begin{aligned}\lb{eq:2.5}
W_{n+}^\eps  f&:= i^n \int_{0\leq t_1 \leq \ldots \leq t_n} e^{i(t_n-t_{n-1})H_0-\eps (t_n-t_{n-1})} V \ldots \\
& e^{i(t_2 - t_1)H_0 - \eps (t_2 - t_1)} V e^{it_1 H_0 - \eps  t_1} V e^{-it_nH_0} f \dd t_1 \ldots \dd t_n,
\end{aligned}\ee
where $\eps>0$.  By \cite[Lemma 4.3]{BeSc},  $W_{n+}^{\eps}\to W_{n+}$  in the strong $L^{2}$ sense as $\eps\to0+$. 
One has the following representation formulas for each $n\ge1$ going back to Yajima, see \cite[Lemma 4.7]{BeSc}: 
\be\begin{aligned}\lb{eq:2.19}
\langle W_{n+}^\eps  f, g \rangle &= (-1)^n \int_{\set R^9} \mc F_{x_{0}}^{-1}T_{n+}^\eps (0, x, y) f(x-y) \ov g(x) \dd y \dd x
\end{aligned}\ee
where for any $\eps  > 0$, $T_{1\pm}^\eps (x_0, x_1, y)$ is defined in the sense of  distributions as
\be\lb{2.16}
(\mc F_{x_0}^{-1} \mc F_{x_1, y} T_{1\pm}^\eps )(\xi_0, \xi_1, \eta) := \frac {\widehat V(\xi_1 - \xi_0)}{|\xi_1 + \eta|^2 - |\eta|^2 \pm i\eps }
\ee
and, more generally, for all $n\ge1$ we have 
\be\begin{aligned}\lb{eq:2.17}
(\mc F_{x_0}^{-1} \mc F_{x_n, y} T_{n\pm}^\eps )(\xi_0, \xi_n, \eta) &:=  
\int_{\R^{3(n-1)}} \frac{\prod_{\ell=1}^n \widehat V(\xi_{\ell} - \xi_{\ell-1}) \dd \xi_1 \ldots \dd \xi_{n-1}}{\prod_{\ell=1}^n (|\xi_{\ell}+\eta|^2-|\eta|^2 \pm  i \eps )}.
\end{aligned}\ee
Even though the first variable does not play a role in~\eqref{eq:2.19}, it is essential in order to express $T_{n+}^{\eps}$ in terms of $T_{1+}^{\eps}$ by means of a convolution structure. In fact,  we  (formally) compose three variable kernels $T(x_0,x_1,y)$ on $\R^9$ by the rule 
\be\begin{aligned}\lb{eq2.21}
(T_1 \oast T_2)(x_0, x_2, y) = \int_{\set R^6} T_1(x_0, x_1, y_1) T_2(x_1, x_2, y-y_1) \dd x_1 \dd y_1.
\end{aligned}\ee
Dually (on the Fourier side), $\oast$ is given by
\EQ{\lb{comp}
&(\mc F^{-1}_{x_0} \mc F_{x_2, y} (T_1 \oast T_2))(\xi_0, \xi_2, \eta) \\
& = \int_{\set R^3} (\mc F^{-1}_{x_0} \mc F_{x_1, y} T_1)(\xi_0, \xi_1, \eta) (\mc F^{-1}_{x_1} \mc F_{x_2, y} T_2)(\xi_1, \xi_2, \eta) \dd \xi_1.
}
So $\oast$ consists of convolution in the $y$ variable  (or multiplication in the dual variable $\eta$),  and composition of operators  relative to the other two variables. 
In the dual coordinates $\xi_0$, $\xi_1$, and $\xi_2$,   composition of operators is preserved.  
Note the order of the variables: $x_0$ is the ``input", $x_2$ the ``output" variable, whereas $y$ is the dual energy variable.

\begin{lemma}\lb{lem:Tast} 
Let $V$ be a Schwartz potential. 
For all $\eps>0$ and  any $n,m\ge1$
$$
T_{m+}^\eps  \oast T_{n+}^\eps  = T_{(m+n)+}^\eps
$$
in the sense of \eqref{comp}. 
\end{lemma}
\begin{proof}
By inspection
\EQ{
\mc F^{-1}_{x_0} \mc F_{x_2, y} T_{2+}^\eps(\xi_0, \xi_2, \eta) &= \int_{\set R^3} \frac {\widehat V(\xi_2-\xi_1)}{|\xi_2 + \eta|^2 - |\eta|^2+i\eps} \cdot \frac {\widehat V(\xi_1-\xi_0)}{|\xi_1 + \eta|^2 - |\eta|^2+i\eps} \dd \xi_1   \\
& = \mc F^{-1}_{x_0} \mc F_{x_2, y} (T_{1+}^\eps\oast T_{1+}^\eps)(\xi_0, \xi_2, \eta)
}
both in the pointwise sense, as well as in the space of distributions.  The general case follows by induction. 
\end{proof}

In order to prove Theorem~\ref{thm:main},   we will show that there exists an algebra under $\oast$ with  the norms of Definition~\ref{def:trip}.  
In~\cite[Section 5]{BeSc} it was shown that 
\EQ{\label{eq:KR}
W_{1+}^\eps  f(x) &=  \int_{\R^3} K_{1+}^\eps (x,x-y) f(y)\, dy \\
K_{1+}^\eps (x,z) &= -\lim_{R\to\infty} \int_{\R^6}  \frac{e^{ix\cdot \xi}\, \what{V}(\xi)\,e^{iz\cdot \eta}}{|\xi+\eta|^2-|\eta|^2+i\eps } e^{-\frac{|\eta|^2}{2R^2}}\,d\xi d\eta
}
Furthermore, 
\EQ{\label{eq:KL}
K_{1+}^\eps(x,z) &=  \const |z|^{-2} e^{-\eps|z|} L(|z|-2x\cdot\hat z, \hat z)
}
where
 for any $\omega\in \Sph^2 $, and $r\in\R$, 
\[
L(r,\omega) = L_{V}(r,\omega):= \int_0^\infty \what V(-s\omega)e^{i\frac{rs}{2}} \,s\, ds
\]
Here $V$ is  any Schwartz function.  The following corollary from~\cite{BeSc} shows how the structure function for $W_{1+}$ arises easily from this formalism.  
It also explains how the norm $\trip\cdot \trip$ arises in Definition~\ref{def:trip}. 

\begin{cor}
\label{cor:K1+}
Let $V$ be a Schwartz function. Define $S_\omega := x-2(\omega\cdot x) \omega$ to be the reflection about the plane $\omega^\perp$. Then for all Schwartz functions $f$ 
\EQ{\label{eq:g1}
(W_{1+}f)(x) = \int_{\Sph^2 }\int_{\R^3} g_1(x,dy,\omega) f(S_\omega x-y)\,  \sigma(d\omega)
}
For fixed $x\in\R^3$, $\omega\in \Sph^2 $ the function $g_1(x,\cdot,\omega)$ is a measure satisfying 
\EQ{\lb{eq:g1meas}
 \int_{\Sph^2 }  \| g_1(x,dy,\omega)\|_{\mes_{y} L^\infty_x}  \, d\omega \le \int_{\Sph^2 }\int_{\R} |L_{V}(r,\omega)|\, drd\omega  = \trip V\trip
}
with $\|\cdot \|_{\mes}$ being the total variation norm for Borel measures. 
\end{cor}
\begin{proof}
From eq.~\eqref{eq:KR} and \eqref{eq:KL}, 
\EQ{
(W_{1+}f)(x) &= \int_0^\infty \int_{\Sph^2 } L(r-2\omega\cdot x,\omega) f(x-r\omega)\, drd\omega \\
&= \int_{\Sph^2 } \int_{\R} \one_{[r>-2\omega\cdot x]} L(r,\omega) f(x-2(\omega\cdot x) \omega -r\omega)\, drd\omega
}
Define
\EQ{\lb{eq:g1def}
g_{1}(x,dy,\omega):= \one_{[(y+2x)\cdot\omega>0]} L(y\cdot\omega,\omega) \, \calH^1_{\ell_\omega}(dy) 
}
Here $\ell_\omega:=\{r\omega\:|\: r\in\R\}$ is the line along $\omega$, and $\calH^1_{\ell_\omega}$ is the $1$--dimensional Hausdorff measure on the line $\ell_\omega$. Then \eqref{eq:g1} holds and 
\EQ{\lb{eq:g1sup}
\| g_1(x,dy,\omega)\|_{L^\infty_x} = | L(y\cdot\omega,\omega)| \, \calH^1_{\ell_\omega}(dy) 
}
which implies \eqref{eq:g1meas}. 
\end{proof}

This result does not explain the origin of the other norm, $\| V\|_{B}$ in Definition~\ref{def:trip}. That norm is needed to bound the higher order structure functions $g_{n}$, $n\ge2$.  The remainder of this paper will be devoted to working out the details of this construction.  To end this section, we recall how \cite{BeSc} fails to reach the scaling-invariant space~$\dot B^{\frac12}$ and we explain how $\|V\|_{B}$ is designed to circumvent the exact difficulty responsible for the loss of $\frac12$ power in Theorem~\ref{thm:struct}. 

First, we point out the connection between $\trip\cdot\trip$ and $\|\cdot\|_{\dot B^{\frac12}}$ as given by \cite[Proposition 6.1]{BeSc}. 

\begin{prop}\lb{lemma2.9} 
Let $L=L_V$ be as above, and $V$  a Schwartz function. Then, with $r\in\R$ and $\omega\in \Sph^2 $, 
\EQ{\label{eq:L2V}
\|L(r,\omega)\|_{L^2_{r, \omega}} \les \|V\|_{L^2}
}
and
\be\lb{est_b}
\begin{aligned}
\|L(r,\omega)\|_{L^1_{r, \omega}} \les \sum_{k \in \Z} 2^{k/2} \|\one_{[2^k, 2^{k+1}]}(|r|) L(r,\omega)\|_{L^2_{r, \omega}} \les \|V\|_{\dot B^{\frac12}}. 
\end{aligned}\ee
Moreover,
for any  $0 < \alpha < 1$,  
\be\lb{2.81}
\sum_{k \in \Z} 2^{\alpha k} \|\one_{[2^k, 2^{k+1}]}(|r|) L(r,\omega)\|_{L^2_{r, \omega}} \les \|V\|_{\dot B^\alpha}.
\ee
\end{prop}

The aforementioned loss of a $\frac12$ power occurred in the following estimate (6.9) from~\cite{BeSc}:
\[
\| v(x) K_{1+}^{\eps}(x,y)\|_{L^{1}_{y}B^{\frac12}_{x}}\les \|v\|_{B^{1}}\|V\|_{B^{\frac12}}
\]
In view of \eqref{eq:KL} this is the same as (in the limit $\eps\to0$)
\[
\int_{S^2} \int_0^\infty \| v(x) L_V(t-2x \cdot \omega, \omega) \|_{B^{\frac12}}\, dtd\omega \les \|v\|_{B^{1}}\|V\|_{B^{\frac12}}
\]
We now show how to avoid this loss by means of the norm~\eqref{eq:Bnorm2}.   

\begin{lemma}
\label{lem:tripB}
For Schwartz functions $v,V$ one has 
\EQ{\label{eq:key}
\int_{S^2} \int_0^\infty \trip v(x) L_V(t-2x \cdot \omega, \omega)\trip \, dt d\omega \les \|v\|_B \trip V\trip.
}
\end{lemma}
\begin{proof}
Writing $\omega^\perp(s)=\omega^\perp + s\omega$, we compute 
\EQ{ 
&\int_{S^2} \int_0^\infty \trip v(x) L_V(t-2x \cdot \omega, \omega)\trip \, dt d\omega  \\
&\le \int_{S^2} \int_{-\infty}^\infty \int_0^\infty \trip \delta_{\omega^\perp(s)} v(x) L_V(t-2x \cdot \omega, \omega)\trip \, dt ds d\omega  \\
&\le \int_{S^2} \int_{-\infty}^\infty \int_{-\infty}^\infty  |L_V(t-2s, \omega)|\, \trip \delta_{\omega^\perp(s)} v(x) \trip \, dt ds d\omega \\
&\le \int_{S^2} \int_{-\infty}^\infty |L_V(t , \omega)|\, dt d\omega \;  \sup_{\omega\in\Sph^2} \int_{-\infty}^\infty \trip \delta_{\omega^\perp(s)} v(x) \trip \,   ds\\
&= \trip V\trip \|v\|_B 
}
which is \eqref{eq:key}. 
\end{proof}

Before continuing with the main argument, the following section exhibits norms that dominate those  in Definition~\ref{def:trip}, but which are more explicit. 
We also check that $\|\cdot\|_{B}$ is scaling invariant. 

\section{A closer look at the norms of Definition~\ref{def:trip}} 
\label{sec:scalinv}

The following lemma bounds the rather implicit $\|\cdot\|_B$-norm by a more explicit Sobolev-type norm. It appears that this cannot be improved significantly.

\begin{lemma}\label{lem:Bnorm dom}
The norm $\|\cdot \|_B$ is scaling invariant  in the sense that with $v_\lambda(x):=\lambda^{-2} v(x/\lambda)$, $\lambda>0$,  one has $\|v\|_B=\|v_\lambda\|_B$. 
Furthermore, 
\EQ{\label{eq:B*}
\|v\|_B\le C\sup_{\omega\in \Sph^2} \int_{-\infty}^\infty \sum_{k\in\Z} 2^{\frac{k}{2}} \big\| \psi(2^{-k} x') v(x'+s\omega) \|_{\dot H^{\frac12}(\omega^\perp)}\, ds =:\|v\|_B^*
}
where $x'\in \omega^\perp$ and $\sum_{k\in\Z} \psi(2^{-k} x') =1 $ for $x'\in \omega^\perp\setminus\{0\}$ is  a Littlewood-Paley partition of unity.  The norm $\|v\|_B^*$ is finite 
on Schwartz functions, and $\|v\|_B^*=\|v_\lambda\|_B^*$ for $\lambda=2^{-\ell}$, $\ell\in\Z$. 
\end{lemma}
\begin{proof}
One has $L_{V}(t,\omega)=\lambda^{-1}L_{V_{\lambda}}(t\lambda^{-1},\omega)$ whence $\trip V_{\lambda}\trip=\trip V\trip$. 
With $\chi$ a standard bump function on the line, 
\EQ{\nn
\|V_{\lambda}\|_{B} &= 
\sup_{\omega\in\Sph^{2}} \lim_{\delta\to0}\int_{-\infty}^{\infty} \Big\|\delta^{-1}\chi((x\cdot\omega-t)/\delta)\lambda^{-2} V(x/\lambda)\Big\|_{L}\, dt\\
&= \sup_{\omega\in\Sph^{2}} \lim_{\delta\to0}\int_{-\infty}^{\infty} \lambda^{-1}\Big\|(\delta/\lambda)^{-1}\chi(((x/\lambda)\cdot\omega-t/\lambda)/(\delta/\lambda))\lambda^{-2} V(x/\lambda)\Big\|_{L}\, dt \\
&= \sup_{\omega\in\Sph^{2}} \lim_{\delta\to0}\int_{-\infty}^{\infty} \Big\|\delta^{-1}\chi(((x/\lambda)\cdot\omega-t)/\delta)\lambda^{-2} V(x/\lambda)\Big\|_{L}\, dt  \\
&=\sup_{\omega\in\Sph^{2}} \lim_{\delta\to0}\int_{-\infty}^{\infty} \Big\|\delta^{-1}\chi((x\cdot\omega-t)/\delta) V(x)\Big\|_{L}\, dt 
= \|V\|_{B}
}
which is the scaling invariance of the $B$-norm. 

To prove \eqref{eq:B*} we fix the plane $\Pi$ to be $x_{1}=0$, or equivalently we set $\omega=(1,0,0)$.  One has, with $v_{s}=\delta_{\Pi(s)}v$, 
\[
L_{v_{s}}(t,\omega)= \int_{0}^{\infty} \mc F_{x_{2},x_{3}} v(s,-\omega_{2}\tau, -\omega_{3}\tau) e^{i\frac{\tau}{2}(t-2s\omega_{1})} \, \tau\,d\tau 
\]
so that 
\EQ{\nn
& \|v\|_{B} = \intR \Big\| \int_{0}^{\infty} \mc F_{x_{2},x_{3}} v(s,-\omega_{2}\tau, -\omega_{3}\tau) e^{i\frac{\tau}{2} t} \, \tau\,d\tau  \Big\|_{L^{1}_{t,\omega}}\, ds \\
&= \intR \!\!\! ds\intR \!\!\!dt\int_{0}^{\pi}\int_{0}^{2\pi} \Big| \int_{0}^{\infty} \mc F_{x_{2},x_{3}} v(s,-\tau\sin\theta\sin\phi, -\tau\sin\theta\cos\phi) e^{i\frac{\tau}{2} t} \, \tau\,d\tau \Big|\,\sin\theta\,  d\phi d\theta\\
&= {\pi} \intR \!\!\! ds\intR \!\!\!dt\int_{0}^{\pi} \Big| \int_{0}^{\infty} \mc F_{x_{2},x_{3}} v(s,-\tau \sin\phi, -\tau \cos\phi) e^{i\frac{\tau}{2} t} \, \tau\,d\tau \Big|\, d \phi \\
&\le C \sum_{k\in\Z}  \intR  2^{\frac{k}{2}} \Big( \intR \int_{0}^{\pi} \one_{[|t|\simeq 2^{k}]} \Big| \int_{0}^{\infty} \mc F_{x_{2},x_{3}} v(s,-\tau \sin\phi, -\tau \cos\phi) e^{i\frac{\tau}{2} t} \, \tau\,d\tau \Big|^{2}\, d \phi  dt\Big)^{\frac12} \, ds
}
For a Schwartz function $w$ in $\R^{2}$ define the sublinear operator $A_{k}w$ as 
\EQ{
A_{k}w = \Big( \intR \int_{0}^{\pi} \one_{[|t|\simeq 2^{k}]} \Big| \int_{0}^{\infty} \what{w}(\tau \sin\phi, \tau \cos\phi) e^{i\frac{\tau}{2} t} \, \tau\,d\tau \Big|^{2}\, d \phi  dt\Big)^{\frac12}
}
Then, on the one hand,
\EQ{
\sum_{k\in\Z}(A_{k}w)^{2}  &\les \intR \int_{0}^{\pi}  \Big| \int_{0}^{\infty} \what{w}(\tau \sin\phi, \tau \cos\phi) e^{i\frac{\tau}{2} t} \, \tau\,d\tau \Big|^{2}\, d \phi  dt
\\
& \les \int_{0}^{\pi} \int_{0}^{\infty} |\what{w}(\tau \sin\phi, \tau \cos\phi)|^{2} | \tau|^{2}\,d\tau d\phi = \| w\|_{\dot H^{\frac12}}^{2}
}
and, on the other hand, 
\begin{align}\nn
&\sum_{k\in\Z} 2^{2k} (A_{k}w)^{2}  \les \intR \int_{0}^{\pi}  \Big| t\int_{0}^{\infty} \what{w}(\tau \sin\phi, \tau \cos\phi) e^{i\frac{\tau}{2} t} \, \tau\,d\tau \Big|^{2}\, d \phi  dt
\\
& \les \int_{0}^{\pi} \int_{0}^{\infty} |\partial_{\tau}(\tau \what{w}(\tau \sin\phi, \tau \cos\phi))|^{2} \,d\tau d\phi  \nn \\
&\les \int_{0}^{\pi} \int_{0}^{\infty} \big[ \tau^{-1} |\what{w}(\tau \sin\phi, \tau \cos\phi))|^{2} + |(\partial_{r}\what{w})(\tau \sin\phi, \tau \cos\phi)|^{2}  \tau\big]\tau\,d\tau d\phi \nn \\
&\les \| \what{w}\|_{\dot H^{\frac12}}^{2} + \| |\xi|^{\frac12}\partial_{r}\what{w}\|_{2}^{2}  \label{eq:xhalf}
\end{align}
The first term in the last line is obtained by Hardy's inequality in the $\xi$ variable, and we bound it further by applying Hardy's inequality in the $x$ variable:
\EQ{\label{eq:x1}
\| \what{w}\|_{\dot H^{\frac12}} = \| |x|^{\frac12} w\|_{2} = \| |x|^{-\frac12} |x|w\|_{2}\les \| |x|w\|_{\dot H^{\frac12}}  
}
For second term in \eqref{eq:xhalf}  we first rewrite $\partial_{r}\what{w}$ as
\EQ{
\partial_{r}\what{w}(\xi) &=  \frac{\xi}{|\xi|} \nabla_{\xi} \int_{\R^{2}} e^{-ix\cdot\xi} w(x)\, dx \\
&= \int_{\R^{2}} e^{-ix\cdot\xi} a(x,\xi) |x| w(x)\, dx,\qquad a(x,\xi)=\frac{x}{|x|}\cdot\frac{\xi}{|\xi|}=:\hat{x}\cdot\hat{\xi}
}
Therefore,
\EQ{
\| |\xi|^{\frac12}\partial_{r}\what{w}\|_{2} &\le \| \hat{x} |x| w\|_{\dot H^{\frac12}}
}
By Lemma~\ref{lem:h12} below one has $\| \hat{x} |x| w\|_{\dot H^{\frac12}}\les \|  |x| w\|_{\dot H^{\frac12}}$. Combining this bound with~\eqref{eq:x1} we conclude that
\EQ{
\sum_{k\in\Z} 2^{2k} (A_{k}w)^{2} \les \|  |x| w\|_{\dot H^{\frac12}}^{2}
}
By  interpolation 
\EQ{
\sum_{k\in\Z} 2^{\frac{k}{2}} A_{k}w \les \| w\|_{*}
}
where in the notation of the real interpolation method
\EQ{
\| w\|_{*} =  \big( \| w\|_{\dot H^{\frac12}}, \| |x| w\|_{\dot H^{\frac12}} \big)_{(\frac12,1)}
}
By Lemma~\ref{lem:h12} the right-hand side is bounded by $\|w\|_{B}^{*}$ and \eqref{eq:B*} is proved.
The other stated properties of $\|\cdot\|_{B}^{*}$ are immediate.  
\end{proof}

The previous proof required two technical properties which we now establish. They are  special cases
of more general statements, but we limit ourselves to what is needed here. 

\begin{lemma}\label{lem:h12}
The following two properties hold:
\begin{itemize}
\item
For any Schwartz function $f$ in $\R^{2}$ one has $\| \hat{x} f\|_{\dot H^{\frac12}}\les \|f\|_{\dot H^{\frac12}}$ where $\hat{x}=x/|x|$. 
\item  With $\|\cdot\|_{B}^{*}$ defined 
as, cf.~\eqref{eq:B*} 
\[
\sum_{k\in\Z} 2^{\frac{k}{2}} \big\| \psi(2^{-k} x) w(x) \|_{\dot H^{\frac12}(\R^{2})}=: \|w\|_{B}^{*} 
\]
one has 
\[
\big( \| w\|_{\dot H^{\frac12}}, \| |x| w\|_{\dot H^{\frac12}} \big)_{(\frac12,1)} \les \|w\|_{B}^{*}
\]
\end{itemize}
\end{lemma}
\begin{proof}
The first property cannot simply be obtained by interpolating between the obvious $L^{2}$ property and the corresponding $\dot H^{1}$ inequality. Indeed, the latter would require Hardy's inequality in $\R^{2}$ with an $r^{-1}$ weight, which fails. So we proceed differently.  Using polar coordinates and complex notation we expand $f$ into a Fourier series:
\[
f(re^{i\theta}) = \sum_{n\in\Z}a_{n}(r) e^{2\pi in\theta},\quad a_{n}(r)= \int_{0}^{1} f(re^{2\pi i\theta}) e^{-2\pi in\theta}\,d\theta
\]
By Plancherel
\EQ{
\|f\|_{2}^{2} = \const \sum_{n\in\Z} \int_{0}^{\infty} |a_{n}(r)|^{2}r\, dr
}
and 
\EQ{
\| f\|_{\dot H^{1}}^{2} &= \| \partial_{r}f\|_{2}^{2}+ \| r^{-1} \partial_{\theta} f\|_{2}^{2} \\
&= \const \sum_{n\in\Z} \int_{0}^{\infty} \big( |a_{n}'(r)|^{2} + \frac{n^{2}}{r^{2}} |a_{n}(r)|^{2}\big)r\, dr
}
By interpolation, 
\EQ{
\| f\|_{\dot H^{\frac12}}^{2} &\simeq   \sum_{n\in\Z} \int_{0}^{\infty} \big( |(-\partial_{r}^{2})^{\frac14}a_{n}(r)|^{2} + \frac{|n|}{r} |a_{n}(r)|^{2}\big)r\, dr
}
Since $\hat{x}=e^{2\pi i\theta}=:e(\theta)$ we conclude that 
\EQ{
\| e(\theta) f\|_{\dot H^{\frac12}}^{2} &\les   \sum_{n\in\Z} \int_{0}^{\infty} \big( |(-\partial_{r}^{2})^{\frac14}a_{n}(r)|^{2} + \frac{|n+1|}{r} |a_{n}(r)|^{2}\big)r\, dr\\
&\les \| f\|_{\dot H^{\frac12}}^{2} + \int_{0}^{\infty} \frac{|a_{0}(r)|^{2}}{r} r\,dr\label{eq:2t}
}
Since 
\[
a_{0}(r) = \int_{0}^{1} f(re(\theta))\, d\theta, 
\]
the final term in \eqref{eq:2t} is 
\EQ{
\int_{0}^{\infty} \frac{|a_{0}(r)|^{2}}{r} r\,dr &\le \int_{0}^{\infty} \int_{0}^{1} \frac{|f(re(\theta))|^{2}}{r} r\,d\theta dr   \\
&\les  \| r^{-\frac12} f\|_{2}^{2} \les \|f\|_{\dot H^{\frac12}}^{2}
}
by Hardy, and the first claim is proved. 

To prove  the second claim we first dominate the weighted norm via a smooth Littlewood-Paley partition of unity, viz.
\EQ{
 \| |x| w\|_{\dot H^{\frac12}} &\les \sum_{k\in\Z} \|\psi(2^{-k}x)|x| w\|_{\dot H^{\frac12}} \les \sum_{k\in\Z} 2^{k} \|\psi(2^{-k}x)  w\|_{\dot H^{\frac12}}
}
For the final inequality it suffices to verify the case $k=0$ by scaling. Then,  by the fractional Leibnitz rule and with $\wtild \psi\psi=\psi$ another Littlewood-Paley function, 
\EQ{
\|\psi(x)|x| w\|_{\dot H^{\frac12}} &\les \| |x|\wtild\psi(x)\|_{\infty}\|\psi(x) w\|_{\dot H^{\frac12}} +  \| |\nabla|^{\frac12 }|x|\wtild\psi(x)\|_{L^{4}}\|\psi(x) w\|_{L^{4}}\\
&\les  \|\psi(x) w\|_{\dot H^{\frac12}}
}
where the final step is obtained by Sobolev embedding.  Clearly, 
\EQ{
 \| w\|_{\dot H^{\frac12}} & \les \sum_{k\in\Z}  \|\psi(2^{-k}x)  w\|_{\dot H^{\frac12}}
}
By the real interpolation property, see \cite[Section 2]{BeSc}, \cite{bergh}
\EQ{
\big( \| w\|_{\dot H^{\frac12}}, \| |x| w\|_{\dot H^{\frac12}} \big)_{(\frac12,1)} &\les \Big( \sum_{k\in\Z}  \|\psi(2^{-k}x)  w\|_{\dot H^{\frac12}}, 
\sum_{k\in\Z} 2^{k} \|\psi(2^{-k}x)  w\|_{\dot H^{\frac12}}\Big)_{(\frac12,1)} \\
&\simeq \sum_{k\in\Z} 2^{\frac{k}{2}} \|\psi(2^{-k}x)  w\|_{\dot H^{\frac12}} = \|w\|_{B}^{*}
}
and we are done. 
\end{proof}

\section{The convolution algebra and the proof of Theorem~\ref{thm:main}}

We now present the algebra formalism in the  scaling invariant setting.   

\begin{definition}
\label{def:Z}
The Banach space $Z$  of tempered distributions is defined as 
\EQ{
Z := \{T(x_0, x_1, y)\in \calS'(\R^{9}) \mid \;&   \calF_{y} T(x_{0},x_{1},\eta) \in L^{\infty}_{\eta}L^{\infty}_{x_{1}}L^{1}_{x_{0}} \} 
}
with norm 
\EQ{\lb{eq:Znorm}
\| T\|_Z := \sup_{\eta\in\R^3} \| \mc F_{y}T(x_0, x_1, \eta)\|_{{L}^{\infty}_{x_1}L^{1}_{x_{0}}} 
}
  $\sup$ being  the essential supremum.   We add the identity $I$ to $Z$, which corresponds to the kernel $T=\delta_0(y)\delta_0(x_1-x_0)$.   The convolution $\oast$ on $T_{1},T_{2}\in Z$ is defined by 
\EQ{\label{eq:FToast}
(T_{1}\oast T_{2} )(x_0, x_2,y) = \mc F_{\eta}^{-1} \Big[ \int_{\R^{3}}\mc F_{y}T_{1}(x_0, x_1, \eta) \mc F_{y}T_{2}(x_1, x_2, \eta)\, dx_{1} \Big](y)
}
\end{definition}

\begin{lemma}\lb{lm2.4} Let
$Z$ is a Banach algebra under $\oast$ with identity element~$I$.
If $V \in L^{3/2, 1}$ then $T_{1+}^\epsilon$ defined by (\ref{2.16}) belongs to $Z$ and  
$\mc F_y T_{1+}^\epsilon$ is given by
\be\lb{2.90}\begin{aligned}
\mc F_y T_{1+}^\epsilon(x_0, x_1, \eta)  
&=  e^{-ix_1 \eta}\,  R_0(|\eta|^2-i\epsilon)(x_0, x_1)V(x_0)\,  e^{ix_0 \eta}.
\end{aligned}\ee
Moreover, 
\EQ{\label{eq:T1Zbd}
\sup_{{\eps>0}}\|T_{1+}^\eps\|_{Z}\les \|V\|_{L^{3/2, 1}} \les \|V\|_{\dot B^{\frac12}}
}
If, in addition, $\|V\|_{L^{3/2, 1}}$ is sufficiently small, then $T_+^\epsilon$ also belongs to $Z$ and
\be\lb{2.37}
(I + T_{1+}^\epsilon) \oast (I - T_+^\epsilon) = (I - T_+^\epsilon) \oast (I + T_{1+}^\epsilon) = I.
\ee
\end{lemma}
\begin{proof}
$Z$ is  a Banach space. The expressions in \eqref{eq:FToast} appearing in brackets  satisfies 
\EQ{
&\sup_{\eta\in\R^{3}}\Big \|  \int_{\R^{3}}\mc F_{y}T_{1}(x_0, x_1, \eta) \mc F_{y}T_{2}(x_1, x_2, \eta)\, dx_{1}  \Big\|_{L^{\infty}_{x_{2}}L^{1}_{x_{0}}} \\
&\le \| \mc F_{y}T_{1}\|_{{L}^{\infty}_{\eta}{L}^{\infty}_{x_1}L^{1}_{x_{0}}} \| \mc F_{y}T_{2}\|_{{L}^{\infty}_{\eta}{L}^{\infty}_{x_2}L^{1}_{x_{1}}} 
= \| T_{1}\|_{Z} \|T_{2}\|_{Z}
}
and so it is a tempered distribution in $\R^{9}$. Therefore, the composition~\eqref{eq:FToast} is well-defined in $Z$ and 
\EQ{\nn
\| T_{1}\oast T_{2}\|_{Z}\le \| T_{1}\|_{Z} \|T_{2}\|_{Z}
}
whence $Z$ is a Banach algebra under $\| \cdot\|_{Z}$.

Formula~\eqref{2.90} follows from \eqref{2.16} by taking Fourier transforms.

By the resolvent identity 
\be\lb{2.29}
(I + R_0(z) V)^{-1} = I - R_V(z) V;\ R_V(z) = (I + R_0(z) V)^{-1} R_0(z)
\ee
for $\Im z>0$. Here $R_{V}(z)=(H-z)^{-1}$ which exists for $\Im z> 0$ since $H$ is self-adjoint.  For $V$ small in $L^{\frac32,1}$ it follows that 
\[
\sup_{\Im z>0}\| R_0(z) Vf\| \le C\|V\|_{L^{\frac32,1}} \|f\|_{\infty} \le \frac12 \|f\|_{\infty}
\]
Hence $(I + R_0(z) V)^{-1}$ exists as a bounded operator on $L^{\infty}$ uniformly in $\Im z>0$, and we may also take the limit $\Im z\ge0$.  In particular,
\[
R_{V}(z) : L^{\frac32, 1}(\R^{3}) \to L^{\infty}(\R^{3})
\]
From (\ref{2.29}), 
\EQ{\label{eq:RI2}
 R_0( |\eta|^2-i\epsilon)V - R_V( |\eta|^2-i\epsilon)V + R_0( |\eta|^2-i\epsilon)V R_V( |\eta|^2-i\epsilon)V  =0
}
whence, with $e^{ix \eta}f(x)=: (M_{\eta}f)(x)$,   
\begin{align}\label{eq:wichtig}
& M_{\eta}^{-1}\, R_0( |\eta|^2-i\epsilon)(x_0,x_1)V(x_0)\,  M_{\eta} - M_{\eta}^{-1} \,  R_V( |\eta|^2-i\epsilon)(x_0,x_1)V(x_0)\,  M_{\eta}\\
&+ M_{\eta}^{-1}\,  R_0( |\eta|^2-i\epsilon)(x_2,x_1)V(x_2)\,  M_{\eta}  \circ M_{\eta}^{-1}\, R_V( |\eta|^2-i\epsilon)(x_0,x_2)V(x_0)\,  M_{\eta}  =0\nn
\end{align}
where $\circ$ signifies integration.  In view of \eqref{2.90} this is tantamount to
\EQ{
0= T_{1+}^\eps - T_{+}^\eps + T_{1+}^\eps\oast T_{+}^\eps
}
or $(I + T_{1+}^\epsilon) \oast (I - T_+^\epsilon)=I$. The second identity in \eqref{2.37} is valid since the resolvent identity also implies \eqref{eq:RI2} with $R_0$ and $R_V$ reversed:
\EQ{\label{eq:RI3}
 R_0( |\eta|^2-i\epsilon)V - R_V( |\eta|^2-i\epsilon)V + R_V( |\eta|^2-i\epsilon)V R_0( |\eta|^2-i\epsilon)V  =0
}
and so that same argument as before concludes the proof. 
\end{proof}

The following spaces play a key role in the proof of Theorem~\ref{thm:main}. The $Y$-space in particular allows us to inductively bound the structure function of each $W_{n+}^\eps$. 
 
\begin{definition}\label{def:VB*}
Let $\cB$ be the closure of the Schwartz functions in $\R^{3}$ under the norm $\trip\cdot\trip$.  
Fix any measurable function $v:\R^3\to \R$ which does not vanish a.e., and so that $\|v\|_{B}<\infty$.   We  introduce the following structures depending on~$v$: 
\begin{itemize}
\item the
 seminormed space
$$
v^{-1}\cB=\{f\ \text{measurable} \mid v(x) f(x) \in \cB \}
$$
with the seminorm $\|f\|_{v^{-1} \cB}:=\trip v f\trip$.  
\item 
Let $X$ the space of two-variable kernels  
\be\begin{aligned}\lb{2.9*}
X := &\Big\{\frak X \in \B(L^{\infty},L^{\infty})\mid (\frak X f)(x) = \int_{\R^3} \frak X(x, y) f(x-y) \dd y,\;\;f\in L^\infty \\
&\|\frak X(x, y)\|_{L^\infty_x L^1_y} < \infty,\ \|\frak X(x, y)\|_{L^1_y v^{-1} \cB_x} < \infty\Big\},
\end{aligned}\ee
with norm (the first $v$ factor is only for homogeneity)
\EQ{\lb{eq:Xnorm*}
\|\frak X\|_X &:= \|v\|_{B} \|\frak X\|_{L^\infty_x L^1_y} + \int_{\R^3} \|v(x) \frak X(x, y)\|_{\cB_x } \dd y.
}
\item 
Let $Y$ be the space of three-variable kernels
\be\lb{2.31*}\begin{aligned}
Y &:= \Big\{T(x_0, x_1, y) \in Z \mid \; \forall f \in  L^\infty \cap v^{-1} \cB\\
& (fT)(x, y):= \int_{\R^3} f(x_0) T(x_0, x_1, y) \dd x_0 \in X_{x_1,y}\Big\},
\end{aligned}\ee
with norm
\EQ{\lb{eq:Ynorm*}
\|T\|_Y &:=  \|T\|_Z  + \|T\|_{B(v^{-1}\cB_{x_0}, X_{x_1, y})}   
}
We   adjoin  an identity element to $Y$, in the form of
\be\lb{ident*}
I(x_0, x_1, y) = \delta_{x_0}(x_1) \delta_0(y) = \delta_{x_1}(x_0) \delta_0(y).
\ee
\end{itemize}
\end{definition} 

Notice that in \eqref{eq:Xnorm*} we use the stronger norm $\|v\|_{B}$ rather than $\trip v\trip$. The presence of 
 $\|T\|_{Z}$ in \eqref{eq:Ynorm*} will require us to ensure that $\|V\|_{B}<\infty$ as well as $V\in\dot B^{\frac12}$. 
 
\begin{lemma}
\label{lem:K1*}
Let $V$ be a Schwartz function, and let $K_{1+}^{\eps}$ be defined in terms of $V$. Then uniformly in $\eps>0$, 
\begin{align}
\| K_{1+}^\eps (x,y)\|_{L^\infty_x L^1_y} & \les \trip V\trip  \lb{eq:Kest1*}  \\
\| v(x) K_{1+}^\eps (x,y)\|_{ L^1_y \cB_x} & \les \| v\|_{B} \trip V\trip     \lb{eq:Kest2*} 
\end{align}
for any $v\in B^1$.  With $f$ a Schwartz function, define a kernel 
\EQ{
\wtil K_{1+}^\eps (x,y) &=  \int_{\R^3} f(x_0)   T_{1+}^\eps (x_0,x,y)\, dx_0
}
with the integral being understood as distributional duality pairing. 
Then  uniformly in $\eps>0$,
\begin{align}
\| \wtil K_{1+}^\eps (x,y)\|_{L^\infty_x L^1_y} & \les \trip fV\trip \lb{eq:Kest1f*}  \\
\| v(x) \wtil K_{1+}^\eps (x,y)\|_{ L^1_y \cB_x} & \les \| v\|_{B} \trip fV\trip   \lb{eq:Kest2f*} 
\end{align}
for any $v\in B^1$. 
\end{lemma}
\begin{proof}
From \eqref{eq:KL} one has for all $\eps>0$, 
\EQ{
\| K_{1+}^\eps (x,y)\|_{L^\infty_x L^1_y} &\le  \const \int_{-\infty}^\infty \int_{\Sph^2} |L(t,\omega)|\, dtd\omega =\const \trip V\trip 
}
which is \eqref{eq:Kest1*}.  For the second estimate \eqref{eq:Kest2*} we invoke Lemma~\ref{lem:tripB}, viz.  
\begin{align}
  \| v(x) K_{1+}^\eps (x,y)\|_{ L^1_y \cB_x} & = \const \int_{\R^3}  \trip v(x) |u|^{-2} L(|u|-2\hat u \cdot x, \hat u)   \trip  \, du \nn \\ 
&= \const \int_{\Sph^2} \int_0^\infty  \trip v(x)   L(t-2\omega \cdot x, \omega)    \trip  \, dt d\omega   \\
&\les  \|v\|_{B}\trip V\trip
\end{align}
as claimed. 

Next, 
\EQ{\label{eq:fV1}
\mc F_{x,y}\wtil K_{1+}^\eps (\xi_1,\eta) &=  \int_{\R^3} \hat{f}(\xi_0) \mc F^{-1}_{x_0} \mc F_{x,y}T_{1+}^\eps (\xi_0,\xi_1,\eta)\, d\xi_0 \\
&= \int_{\R^3} \frac{ \hat{f}(\xi_0) \what{V}(\xi_1-\xi_0)}{|\xi_1+\eta|^2-|\eta|^2+i\eps}\, d\xi_0 = \frac{ \what{fV}(\xi_1)}{|\xi_1+\eta|^2-|\eta|^2+i\eps}
}
In view of \eqref{2.16}, this leads to  the kernel $K_{1+}^\eps$ associated with the potential $fV$.
\end{proof} 

Next, we define the operation of  {\em contraction}: 

\begin{lemma}
For $\frak X \in X$,  the \emph{contraction} of $T\in Y$ by $\frak X$ is 
\be\lb{contractie}
(\frak X T)(x, y) := \int_{\R^6} \frak X(x_0, y_0) T(x_0, x, y-y_0) \dd x_0 \dd y_0.   
\ee
Then $\frak X T\in X_{x,y}$ and $\|\frak X T\|_{X}\le \|T\|_{Y}\|\frak X\|_{X}$.  We interpret the right-hand side of \eqref{contractie} relative to the Fourier variable: 
\EQ{\lb{eq:mcFXT}
\mc F^{-1}_{\eta}\Big [ \int_{\R^3} \mc F_{y_{0}} \frak X(x_0, \eta) \mc F_{y_{0}}T(x_0, x,\eta) \dd x_0 \Big](y)
}
The integral is absolutely convergent and the inverse Fourier transform relative to $\eta$ is a tempered distribution.  
\end{lemma}
\begin{proof}
We have $\mc F_{y_{0}} \frak X(x_0, \eta) \in L^{\infty}_{x_{0},\eta}$ and $\mc F_{y_{0}}T(x_0, x,\eta) \in L^{\infty}_{x}L^{1}_{x_{0}}$, whence the claim about the 
integral in brackets.  The estimate $\|\frak X T\|_{X}\le \|T\|_{Y}\|\frak X\|_{X}$  follows from the definition of the space $Y$:
\EQ{\nn
& \Big \| \int_{\R^6} \frak X(x_0, y_0) T(x_0, x, y-y_0) \dd x_0 \dd y_0 \Big\|_{X_{x,y}} \\
&\le  \int_{\R^{3}}  \|  (\frak X(\cdot, y_0) T)(x,y-y_{0})  \|_{X_{x,y}}   \dd y_{0}\\
& = \int_{\R^{3}}  \|  (\frak X(\cdot, y_0) T)(x,y)  \|_{X_{x,y}}   \dd y_{0}\\ 
&\les  \|T\|_{Y} \int_{\R^{3}} \| \frak X(\cdot, y_0) \|_{V^{-1}\cB_{x}}   \dd y_{0} \le \|T\|_{Y}\|\frak X\|_{X}
}
and we are done.
\end{proof}

The previous lemma allows us to prove that $Y$ is a Banach algebra under the composition~$\oast$.  This will allow us to prove the key property 
that $T_{n+}^{\eps}\in Y$ starting from the case $T_{1+}\in Y$,  which we now state. 

\begin{lemma}\lb{yalg}
$Y$ defined by (\ref{2.31*}) is a  Banach algebra with the operation $\oast$ defined in the ambient algebra~$Z$.  
\end{lemma}
\begin{proof}
The fact that $\oast$ is associative (and non-commutative) is clear in $Z$, and the unit element is given by (\ref{ident*}). Since $Y \subset Z$, the same is true in $Y$.

 The definitions of $X$ and $Y$ imply that each contraction $\frak X T$ (see (\ref{contractie})) is in $X$ and $\|\frak X T\|_X \les \|\frak X\|_X \|T\|_Y$.
 We have 
\EQ{\label{eq:T3oast}
&\int_{\R^3} f(x_0) T_3(x_0, x_2, y) \dd x_0 = \\
&= \int_{\set R^9} f(x_0) T_1(x_0, x_1, y_1) T_2(x_1, x_2, y-y_1) \dd x_1 \dd y_1 \dd x_0.
}
As in the case of \eqref{contractie}, the $y$-integral is to be understood in the distributional Fourier sense. 
Integrating in $x_0$, we obtain an expression of the form $\frak X T_2$ for $\frak X \in X$ with $\|\frak X\|_X \les \|f\|_{V^{-1} B} \|T_1\|_Y$. Then $\frak X T_2$ belongs to $X$ as stated above and has a norm  at most $\les \|f\|_{V^{-1} B} \|T_1\|_Y \|T_2\|_Y$. Thus, $T_3 = T_1 \oast T_2 \in Y$ and
$$
\|T_1 \oast T_2\|_Y \le C \|T_1\|_Y \|T_2\|_Y
$$
with some absolute constant $C$. Multiplying the norm by $C$ removes this constant from the previous inequality, and so $Y$ is an algebra under this new norm. 
\end{proof}

\begin{corollary}
\label{cor:T1Y*}
Let $V$ be Schwartz and apply Definition~\ref{def:VB*} with $v=V$, the potential. Then for every $\eps>0$ we have $T_{1+}^{\eps}\in Y$  and
\EQ{
\sup_{\eps>0} \| T_{1+}^{\eps} \|_{Y} \les \| V\|_{B}+ \|V\|_{\dot B^{\frac12}}
}
\end{corollary}
\begin{proof}
By \eqref{eq:T1Zbd} we have 
\[
\sup_{\eps>0} \| T_{1+}^{\eps}\|_{Z}\les \|V\|_{\dot B^{\frac12}}
\]
It remains to show that 
\EQ{
\sup_{\eps>0}\Big\| \int_{\R^{3}} f(x_{0}) T_{1+}^{\eps}(x_{0},x,y)\, dx_{0}\Big\|_{X_{x,y}} \les  \|V\|_{B} \trip fV\trip   
}
In view of \eqref{eq:Xnorm*} this is implied by Lemma~\ref{lem:K1*}. 
\end{proof}

We are now in a position to obtain the key representation result concerning the partial wave operators $W_{n+}$, see~\eqref{Dyson}. 
In what follows, we let $B_{*}$ be the space obtained as the closure of the Schwartz functions under the norm 
$$\|\cdot\|_{B_{*}}:=\|\cdot\|_{B} +\|\cdot\|_{\dot B^{\frac12}},$$
see \eqref{eq:B*}.  We define $B_*$ as the space obtained as the closure of Schwartz functions under the norm $\|\cdot\|_{B_{*}}$. 

\begin{prop}
Let $V$ be a Schwartz potential. Then $T_{n+}^{\eps}\in Y$ for any $n\ge1$ and $\eps>0$ and 
\EQ{\label{eq:WnY*}
\sup_{\eps>0} \| T_{n+}^{\eps}\|_{Y}  \le C^{n}\|V\|_{B_{*}}^{n}
}
with some absolute constant $C$.  Moreover, for all Schwartz functions $f$ one has 
\EQ{\label{eq:gn*}
(W_{n+}^{\eps}f)(x) = \int_{\Sph^2 }\int_{\R^3} g_n^{\eps}(x,dy,\omega) f(S_\omega x-y)\,  \sigma(d\omega)
}
where  for fixed $x\in\R^3$, $\omega\in \Sph^2 $ the expression $g_n^{\eps}(x,\cdot,\omega)$ is a measure satisfying 
\EQ{\lb{eq:gnmeas*}
\sup_{\eps>0} \int_{\Sph^2 }  \| g_n^{\eps}(x,dy,\omega)\|_{\mes_{y} L^\infty_x}  \, d\omega \le  C^{n}\|V\|_{B_{*}}^{n}
}
where $\|\cdot \|_{\mes}$ refers to the total variation norm of Borel measures.   The same conclusion also holds if $V\in B_*$. 
\end{prop}
\begin{proof}
First,  $T_{n+}^{\eps} = T_{1+}^{\eps}\oast T_{(n-1)+}^{\eps}$. Corollary~\ref{cor:T1Y*} and the algebra property of $Y$ imply~\eqref{eq:WnY*} by induction. 
Second, we have
\EQ{\lb{eq:1Tn*}
W_{n+}^{\eps} &= (-1)^{n} \one_{\R^{3}}T_{n+}^{\eps} = (-1)^{n} \one_{\R^{3}}(T_{(n-1)+}^{\eps}\oast T_{1+}^{\eps} )\\
&= -((-1)^{n-1} \one_{\R^{3}} T_{(n-1)+}^{\eps}) T_{1+}^{\eps} = - W_{(n-1)+}^{\eps} T_{1+}^{\eps}
}
The notation in the second line contraction of a kernel in $Y$ by an element of $X$; this follows again by induction starting from  $W_{0+}^{\eps}=\one_{\R^{3}}$ via~\eqref{contractie}.  By the boundedness of $T_{n+}^{\eps}$ in $Y$ it follows that the right-hand side of \eqref{eq:1Tn*} is well-defined in $Y$.  Thus,  by the first equality sign in~\eqref{eq:1Tn*}, 
\EQ{\label{eq:WnX*}
\sup_{\eps>0} \| W_{n+}^{\eps}\|_{X} &\le \|\one_{\R^{3}}\|_{V^{-1}\cB} \sup_{\eps>0} \|T_{n+}^{\eps}\|_{Y} \le C^{n+1}\trip V\trip \|V\|_{B_{*}}^{n}\\
&\le C^{n+1} \|V\|_{B_{*}}^{n+1}
}
We denote the kernel of $W_{1+}^{\eps}$ by $\frak X_{V}^{\eps}$, where $V$ is the potential.  Thus,
\[
\frak X_{V}^{\eps} (x,y) = - \int_{\R^{3}}  T_{1+}^{\eps} (x_{0},x,y)\, dx_{0} =  -(\one_{\R^{3}} T_{1+}^{\eps})(x,y) \in X
\]
By  \eqref{eq:1Tn*}, 
\EQ{
W_{n+}^{\eps}(x,y) &= -\int_{\R^{6}} W^{\eps}_{(n-1)+}(x',y')T_{1+}^{\eps}(x',x,y-y') \, dx'dy' \\
&= \int_{\R^{3}} \frak X^{\eps}_{f^{\eps}_{y'}V}(x,y-y')\, dy'
}
Here we wrote $f^{\eps}_{y'}(x')=W^{\eps}_{(n-1)+}(x',y')$ and we used \eqref{eq:fV1}. 

We now invoke the representation from Corollary~\ref{cor:K1+}. Specifically, by \eqref{eq:g1} there exists $g_{1,f^{\eps}_{y'}}^{\eps}(x,dy,\omega)$  so that for every 
$\phi\in \mc S$ one has
\EQ{\nn
(\frak X^{\eps}_{f^{\eps}_{y'}V}\; \phi)(x) = \int_{\Sph^2 }\int_{\R^3} g_{1,f^{\eps}_{y'}}^{\eps} (x,dy,\omega) \phi(S_\omega x-y)\,  \sigma(d\omega)
}
where  for fixed $x\in\R^3$, $\omega\in \Sph^2 $ the expression $g_{1,f^{\eps}_{y'}}^{\eps}(x,\cdot,\omega)$ is a measure satisfying 
\EQ{\nn
\sup_{\eps>0} \int_{\Sph^2 }  \| g_{1,f^{\eps}_{y'}}^{\eps}(x,dy,\omega)\|_{\mes_{y} L^\infty_x}  \, d\omega &\le C\trip f^{\eps}_{y'}V\trip \\
&= C  \| W^{\eps}_{(n-1)+}(x',y') \|_{V^{-1}\cB_{x'}}   
}
Therefore, 
\begin{align} 
(W^{\eps}_{n+}\phi)(x) &= \int_{\R^{3}} W_{n+}^{\eps}(x,y) \phi(x-y)\, dy \nn  \\
&= \int_{\R^{6}} \frak X^{\eps}_{f^{\eps}_{y'}V}(x,y-y') \phi(x-y)\, dy dy'  \nn  \\
&=   \int_{\R^{6}} \frak X^{\eps}_{f^{\eps}_{y'}V}(x,y) \phi(x-y-y')\, dy dy'  \nn  \\
& = \int_{\R^{3}} (\frak X^{\eps}_{f^{\eps}_{y'}V} \phi)(x-y')\,dy' \lb{eq:I} \\
&= \int_{\R^{3}} \int_{\Sph^2 }\int_{\R^3} g_{1,f^{\eps}_{y'}}^{\eps} (x-y',dy,\omega) \phi(S_\omega (x-y')-y)\,  \sigma(d\omega) \,dy' \nn \\
&=  \int_{\Sph^2 }\int_{\R^3} \Big[ \int_{\R^{3}} g_{1,f^{\eps}_{y'}}^{\eps} (x-y',d(y-S_{\omega}y'),\omega) \, dy'\Big] \phi(S_\omega x-y)\,  \sigma(d\omega)   \nn 
\end{align}
The expressions in brackets is the structure function
\EQ{\lb{eq:II}
g_{n}(x,dy, \omega) := \int_{\R^{3}} g_{1,f^{\eps}_{y'}}^{\eps} (x-y',d(y-S_{\omega}y'),\omega) \, dy'
}
In fact, it is a measure in the $y$-coordinate and  
\EQ{\lb{eq:III}
(W^{\eps}_{n+}\phi)(x) &= \int_{\Sph^{2}} \int_{\R^{3}}g_{n}(x,dy, \omega) \phi(S_\omega x-y)\,  \sigma(d\omega) 
}
Moreover, we have the bounds, uniformly in $\eps>0$ 
\EQ{\nn
 & \int_{\Sph^2 }  \| g_n^{\eps}(x,dy,\omega)\|_{\mes_{y} L^\infty_x}  \, d\omega  \\
& = \int_{\Sph^2 }\int_{\R^3}    \big\| g_{1,f^{\eps}_{y'}}^{\eps} (x-y',d(y-S_{\omega}y'),\omega) \big \|_{\mes_{y} L^\infty_x}   \,  dy'  d\omega\\
&= \int_{\Sph^2 }\int_{\R^3}   \big\| g_{1,f^{\eps}_{y'}}^{\eps} (x,dy,\omega) \big \|_{\mes_{y} L^\infty_x}   \,  dy'  d\omega   \\
&\le    C  \int_{\R^3}  \| W^{\eps}_{(n-1)+}(x',y') \|_{V^{-1}\cB_{x'}} \,  dy' \\
&  = C   \| W^{\eps}_{(n-1)+}(x',y') \|_{L^{1}_{y'}V^{-1}\cB_{x'}} \le  C   \| W^{\eps}_{(n-1)+} \|_{X} \le C^{n}\|V\|_{B_{*}}^{n}
}
by \eqref{eq:WnX*}. This concludes the argument 
under the assumption that $f^{\eps}_{y'}(x')$ is a Schwartz function. To remove this assumption,  
we can make $$\| W^{\eps}_{(n-1)+}(x',y')-\tilde f^{\eps}_{y'}(x') \|_{X}$$ arbitrarily small with a 
Schwartz function $\tilde f^{\eps}_{y'}(x')$ in $\R^{6}$. Then the previous calculation  shows that 
\[
\int_{\Sph^2 }  \| g_n^{\eps}(x,dy,\omega)- \tilde g_n^{\eps}(x,dy,\omega)\|_{\mes_{y} L^\infty_x}  \, d\omega 
\]
can be made as small as we wish where $\tilde g_n^{\eps}(x,dy,\omega)$ is the function generated by $\tilde f^{\eps}_{y'}(x')$. Passing to the limit concludes the proof. 

To remove the assumption that $V$ be a Schwartz function, we approximate $V\in B_*$ by Schwartz functions in the norm $\|\cdot\|_{B_*}$. 
We achieve convergence of of the functions $g_{n}$ by means of \eqref{eq:gnmeas*} and of the kernels $W_{n+}^{\eps}$ 
themselves by means of~\eqref{eq:WnX*}. To be specific, denoting by $\widetilde W^{\eps}_{n+}$ and  $\tilde g_{n}$ the quantities corresponding to the potential $\tilde V$,
taking differences yields 
\EQ{\nn 
& \| \widetilde  W^{\eps}_{n+} -  W^{\eps}_{n+} \|_{X} + \int_{\Sph^2 } \| g_n^{\eps}(x,dy,\omega)-\tilde g_n^{\eps}(x,dy,\omega)\|_{\mes_{y} L^\infty_x}  \, d\omega \\
& \le C^{n} \| V-\tilde V\|_{B_*}(\|V\|^{n-1}_{B_*}+ \|\tilde V\|_{B_*}^{n-1})
}
uniformly in $\eps>0$.  
\end{proof}

To prove Theorem~\ref{thm:main} we now simply sum the series $\sum_{n=1}^\infty g_n $ which can be done in view of the previous proposition,  provided $c_0$ is sufficiently small.


\begin{thebibliography}{AbcDef1}

 

\bibitem[Bec]{bec} Beceanu, M. \emph{New estimates for a time-dependent Schr\"{o}dinger equation}, Duke Math.\ J.\ Volume 159, Number 3 (2011), pp.\ 417--477.

 

\bibitem[BeGo]{becgol} Beceanu, M., Goldberg, M. \emph{Schr\"{o}dinger dispersive estimates for a scaling-critical class of potentials}, Comm.\  Math.\ Phys.,  Vol.~314 (2012), Issue 2, pp.\ 471--481.

\bibitem[BecSch]{BeSc} Beceanu, M.,  Schlag, W. \emph{ Structure formulas for wave operators}. Preprint 2016. 

\bibitem[BeL\"o]{bergh} Bergh, J., L\"ofstr\"om, J. \emph{Interpolation Spaces. An Introduction}, Springer-Verlag, 1976.


 \bibitem[Kat]{kato} Kato, T. \emph{Wave operators and similarity for some non-selfadjoint operators}, Math.\ Ann.\ 162 (1965/1966), pp.\ 258--279. 


 


\bibitem[Yaj1]{yajima0} Yajima, K. \emph{The $W^{k, p}$-continuity of wave operators for Schr\"{o}dinger operators}, Proc.\ Japan Acad., 69, Ser.\ A (1993), pp.\ 94--99.

\bibitem[Yaj2]{yajima} Yajima, K. \emph{The $W^{k, p}$-continuity of wave operators for Schr\"{o}dinger operators}, J.\ Math.\ Soc.\ Japan 47 (1995), pp.\ 551--581.


 \end{thebibliography}
\end{document}